\documentclass[12pt]{article}

\usepackage{amssymb}
\usepackage{amsmath}
\usepackage[mathscr]{eucal}
\usepackage{cite}
\usepackage[arrow, matrix, curve]{xy}

\begin{document}

\renewcommand{\citeleft}{{\rm [}}
\renewcommand{\citeright}{{\rm ]}}
\renewcommand{\citepunct}{{\rm,\ }}
\renewcommand{\citemid}{{\rm,\ }}

\newcounter{abschnitt}
\newtheorem{theorem}{Theorem}[abschnitt]
\newtheorem{koro}{Corollary}[abschnitt]
\newtheorem{defi}{Definition}
\newtheorem{satz}[koro]{Theorem}
\newtheorem{lem}[koro]{Lemma}

\newcounter{saveeqn}
\newcommand{\alpheqn}{\setcounter{saveeqn}{\value{abschnitt}}
\renewcommand{\theequation}{\mbox{\arabic{saveeqn}.\arabic{equation}}}}
\newcommand{\reseteqn}{\setcounter{equation}{0}
\renewcommand{\theequation}{\arabic{equation}}}

\newcommand{\up}[1]{\stackrel{\raisebox{-0.07cm}{\tiny$\smile$}}{#1}}
\newcommand{\down}[1]{\stackrel{\raisebox{-0.08cm}{\tiny$\frown$}}{#1}}

\hyphenation{convex} \hyphenation{bodies}

\sloppy

\phantom{a}

\vspace{-1.5cm}

\begin{center} \begin{LARGE} Crofton Measures and Minkowski
Valuations \\[0.7cm] \end{LARGE}

\begin{large} Franz E. Schuster \end{large}
\end{center}

\vspace{-0.8cm}

\begin{quote}
\footnotesize{ \vskip 1truecm\noindent {\bf Abstract.} A
description of continuous rigid motion compatible Minkowski
valuations is established. As an application we present a
Brunn--Minkowski type inequality for intrinsic volumes of these
valuations.}
\end{quote}

\vspace{0.6cm}

\centerline{\large{\bf{ \setcounter{abschnitt}{1}
\arabic{abschnitt}. Introduction}}} \alpheqn

\vspace{0.6cm}

As a generalization of the notion of measure, valuations on
convex bodies (compact convex sets) have always played a central
role in geometry. They were the critical ingredient in Dehn's
solution of Hilbert's third problem and they have since been
intimately tied to the dissection theory of polytopes. The
starting point for many important new results in valuation theory
is Hadwiger's \textbf{\cite{hadwiger51}} remarkable
characterization of the intrinsic volumes as the continuous rigid
motion invariant valuations. For more information on the history
of valuations, see \textbf{\cite{Klain:Rota}} and
\textbf{\cite{McMullen93}}. For some of the more recent results,
see, e.g., \textbf{\cite{Alesker99, Alesker01, Alesker03,
Alesker04, Alesker04b, Bernig08a, Bernig03, fu06, klain99,
ludwig02a, Ludwig:matrix, LR99, centro}}.

In 2001, Alesker \textbf{\cite{Alesker01}} has given a complete
description of continuous translation invariant valuations on
convex bodies thereby confirming, in a much stronger form, a
conjecture by McMullen. Alesker's landmark result, now known as
the Irreducibility Theorem, has subsequently led to the discovery
of several new operations on the space of continuous translation
invariant valuations, illuminating a new and rich algebraic
structure \textbf{\cite{Alesker01, Alesker03, Alesker04,
Alesker04b, Bernig03, bernigfu06}}. In a different line of
research, with a similar far reaching impact on the theory of
valuations, Ludwig \textbf{\cite{ludwig02, Ludwig:Minkowski,
Ludwig06}} first obtained characterizations of convex and
\linebreak star body valued valuations which are compatible with
nondegenerate linear transformations, see also
\textbf{\cite{hablud06}}. Her results revealed the underlying
reason why such basic notions as projection, centroid and
intersection bodies are indeed fundamental in the affine theory
of convex bodies.

In this paper, we apply deep results on translation invariant
real-valued valuations to establish a description of even and
translation invariant convex body valued valuations which
intertwine orthogonal transformations. Our result provides a
significant extension of earlier work by Schneider
\textbf{\cite{schneider74}}, Kiderlen \textbf{\cite{kiderlen05}},
and the author \textbf{\cite{Schu06a}}. As an application, we
obtain a new Brunn--Minkowski type inequality for intrinsic
volumes of these valuations, generalizing results by Lutwak
\textbf{\cite{lutwak93}} and the author \textbf{\cite{Schu06b}}.

\pagebreak

Let $\mathcal{K}^n$ denote the space of convex bodies in
$\mathbb{R}^n$, $n \geq 3$, endowed with the Hausdorff metric. A
convex body $K$ is uniquely determined by its support function
$h(K,u)=\max\{u\cdot x: x \in K\}$, for $u \in S^{n-1}$. For $i
\in \{1, \ldots, n - 1\}$, let $\mathrm{Gr}_{i,n}$ be the
Grassmannian of $i$-dimensional subspaces in $\mathbb{R}^n$. The
$i$th projection function $\mathrm{vol}_i(K|\,\cdot\,)$ of $K \in
\mathcal{K}^n$ is the continuous function on $\mathrm{Gr}_{i,n}$
defined such that $\mathrm{vol}_i(K|E)$, for $E \in
\mathrm{Gr}_{i,n}$, is the $i$-dimensional volume of the
orthogonal projection of $K$ onto $E$.

\vspace{0.3cm}

\noindent {\bf Definition} \emph{A map $\Phi: \mathcal{K}^n
\rightarrow \mathcal{K}^n$ is called a Minkowski valuation if
\[\Phi K + \Phi L = \Phi(K \cup L) + \Phi(K \cap L),  \]
whenever $K \cup L \in \mathcal{K}^n$ and addition on
$\mathcal{K}^n$ is Minkowski addition.}

\vspace{0.3cm}

Important examples of Minkowski valuations are such central
notions as the projection and the difference operator: The {\it
projection body} $\Pi K$ of $K$ is the convex body defined by
\[h(\Pi K,u)=\mathrm{vol}_{n-1}(K|u^{\bot}), \qquad u \in S^{n-1},  \]
where $u^{\bot}$ is the hyperplane orthogonal to $u$. The {\it
difference body} $\mathrm{D}K$ of $K$ can be defined by
\[h(\mathrm{D}K,u)=\mathrm{vol}_{1}(K|u), \qquad u \in S^{n-1}.  \]

First results on a special class of Minkowski valuations were
obtained by Schneider \textbf{\cite{schneider74}} in the 1970s,
but only through the recent seminal work of Ludwi\-g
\textbf{\cite{ludwig02, Ludwig:Minkowski}} classifications of
convex and star body valued valuations have become the focus of
increased attention, see \textbf{\cite{haberl08, hablud06,
kiderlen05, Ludwig06, schnschu, Schu06a}}. For example, Ludwig
\textbf{\cite{Ludwig:Minkowski}} established characterizations of
the projection and the difference operator as unique Minkowski
valuations which are compatible with affine transformations of
$\mathbb{R}^n$.

In this article, we consider continuous translation invariant
Minkowski valuations which are $\mathrm{O}(n)$ equivariant. This
class of operators was studied under additional homogeneity
assumptions first by Schneider \textbf{\cite{schneider74}}, and
more recently by Kiderlen \textbf{\cite{kiderlen05}} and the
author \textbf{\cite{Schu06a}}. A map $\Phi$ from $\mathcal{K}^n$
to $\mathcal{K}^n$ \linebreak (or $\mathbb{R}$) is said to have
{\it degree} $i$ if $\Phi(\lambda K)=\lambda^i\Phi K$ for $K \in
\mathcal{K}^n$ and $\lambda \geq 0$. In the case of degree $i$
Minkowski valuations, Kiderlen \textbf{\cite{kiderlen05}}, for $i
= 1$, and the author \textbf{\cite{Schu06a}}, for $i = n - 1$,
recently obtained representations of these maps by spherical
convolution operators. A description of the intermediate degree
cases $i \in \{2, \ldots, n - 2\}$ remained open (by a result of
McMullen \textbf{\cite{McMullen77}}, only integer degrees $0 \leq
i \leq n$ can occur, cf.\ Section 3).

\pagebreak

As our main result we establish a representation for smooth
translation \linebreak invariant and $\mathrm{O}(n)$ equivariant
Minkowski valuations $\Phi_i$ of degree \linebreak $i \in \{1,
\ldots, n - 1\}$ which are {\it even}, i.e., $\Phi_i(-K)=\Phi_i
K$ for $K \in \mathcal{K}^n$. We show that these maps are
generated by convolution of the projection \linebreak functions
with (invariant) measures on the sphere.

\begin{theorem} \label{thm1} Let $\Phi_i: \mathcal{K}^n \rightarrow
\mathcal{K}^n$ be a smooth translation invariant and
$\mathrm{O}(n)$ equivariant Minkowski valuation of degree $i \in
\{1, \ldots, n - 1\}$. If $\Phi_i$ is even, then there exists an
$\mathrm{O}(i) \times \mathrm{O}(n - i)$ invariant measure $\mu$
on $S^{n-1}$ such that for every $K \in \mathcal{K}^n$,
\begin{equation} \label{main}
h(\Phi_iK,\cdot) = \mathrm{vol}_i(K|\,\cdot\,) \ast \mu.
\end{equation}
\end{theorem}

Theorem \ref{thm1} provides an extension of the previously known
convolution representations of Kiderlen
\textbf{\cite{kiderlen05}} and the author \textbf{\cite{Schu06a}}
(cf.\ Section 5). The convolution in (\ref{main}) is induced from
$\mathrm{O}(n)$ by identifying $S^{n-1}$ and $\mathrm{Gr}_{i,n}$
with the homogeneous spaces $\mathrm{O}(n)/\mathrm{O}(n - 1)$ and
$\mathrm{O}(n)/\mathrm{O}(i) \times \mathrm{O}(n - i)$. The
generating measures for the projection and the difference
operator are Dirac measures (cf.\ Section 4; for additional
examples see Section 5).

The invariant signed measures in Theorem \ref{thm1} are Crofton
measures of associated real-valued valuations (cf.\ the proof of
Theorem \ref{thm1}; for related results see
\textbf{\cite{Bernig07}}). Additional properties of these
measures and uniqueness will be discussed in Section 6.

The notion of smooth translation invariant real-valued valuations
was \linebreak introduced by Alesker in
\textbf{\cite{Alesker03}}. We will extend this definition to
translation invariant Minkowski valuations which intertwine
orthogonal transformations in Section 5.

As a consequence of Theorem \ref{thm1}, we obtain in Section 6 a
stronger result, Theorem \ref{stronger}, describing the class of
smooth translation invariant and $\mathrm{O}(n)$ equivariant even
Minkowski valuations without additional assumption on the degree.
We complement these results with the following:

\begin{theorem} \label{thm2} Every continuous translation invariant and $\mathrm{O}(n)$
equivariant even Minkowski valuation can be approximated
uniformly on compact subsets of $\mathcal{K}^n$ by smooth
translation invariant and $\mathrm{O}(n)$ equi\-variant even
Minkowski valuations.
\end{theorem}

Consequently, the problem of describing continuous translation
invariant and $\mathrm{O}(n)$ equivariant even Minkowski
valuations is reduced to a description of smooth ones which is
provided by Theorem \ref{thm1} and Theorem \ref{stronger}.

For $i \in \{1, \ldots, n - 1\}$, let $V_i(K)$ denote the $i$th
intrinsic volume of $K \in \mathcal{K}^n$ and denote by $\Pi_i K$
the projection body of order $i$ defined by
\[h(\Pi_i K,u)=V_i(K|u^{\bot}), \qquad u \in S^{n-1}.   \]
In \textbf{\cite{lutwak85, lutwak93}} Lutwak obtained an array of
geometric inequalities for the intrinsic volumes of projection
bodies which have been recently generalized in
\textbf{\cite{Schu06b}}. (For important recent related results on
$L_p$ projection bodies, see \textbf{\cite{LYZ2000a, LYZ2000b,
LYZ2004}}). A special case of \textbf{\cite[\textnormal{Theorem
6.2}]{lutwak93}} is the following: If $K, L \in \mathcal{K}^n$
have non-empty interior and $i \in \{2, \ldots, n - 1\}$, then
\begin{equation} \label{piiinequ}
V_{i+1}(\Pi_i(K + L))^{1/i(i+1)} \geq
V_{i+1}(\Pi_iK)^{1/i(i+1)}+V_{i+1}(\Pi_iL)^{1/i(i+1)},
\end{equation}
with equality if and only if $K$ and $L$ are homothetic.

As an application of Theorem \ref{thm1} we obtain a similar
Brunn--Minkowski type inequality for all continuous translation
invariant and $\mathrm{O}(n)$ equivariant even Minkowski
valuations of a given degree.

\begin{theorem} \label{thm3} Let $\Phi_i: \mathcal{K}^n \rightarrow
\mathcal{K}^n$ be a continuous translation invariant and
$\mathrm{O}(n)$ equivariant even Minkowski valuation of degree $i
\in \{2, \ldots, n - 1\}$. \linebreak If $K, L \in \mathcal{K}^n$
have non-empty interior, then
\[V_{i+1}(\Phi_i(K + L))^{1/i(i+1)} \geq V_{i+1}(\Phi_iK)^{1/i(i+1)}+V_{i+1}(\Phi_iL)^{1/i(i+1)}.  \]
If $\Phi_i$ maps convex bodies with non-empty interiors to bodies
with non-empty interiors, then equality holds if and only if $K$
and $L$ are homothetic.
\end{theorem}

Note that Theorem \ref{thm3} provides a significant
generalization of Lutwak's inequality (\ref{piiinequ}) and the
related results in \textbf{\cite{Schu06b}}. We also remark that
the classical Brunn--Minkowski inequalities for intrinsic volumes
are special cases of Theorem \ref{thm3}.

\vspace{1cm}

\centerline{\large{\bf{ \setcounter{abschnitt}{2}
\arabic{abschnitt}. Background material}}}

\vspace{0.7cm} \reseteqn \alpheqn

In the following we recall basic facts about convex bodies and
mixed volumes. For quick reference, we state the geometric
inequalities from the Brunn--Minkowski theory needed in the proof
of Theorem 3. For general reference the reader may wish to
consult the book by Schneider \textbf{\cite{schneider93}}.

A convex body $K \in \mathcal{K}^n$ is uniquely determined by the
values of its support function $h(K,\cdot)$ on $S^{n-1}$. From
the definition of $h(K,\cdot)$, it is easily seen that
$h(\vartheta K,u)=h(K, \vartheta^{-1}u)$ for every $u \in
S^{n-1}$ and every $\vartheta \in \mathrm{O}(n)$.

\pagebreak

The {\it Steiner point} $s(K)$ of $K \in \mathcal{K}^n$ is the
point in $K$ defined by
\[s(K)=n\int_{S^{n-1}}h(K,u)u\,du,  \]
where the integration is with respect to the rotation invariant
probability measure on $S^{n-1}$.

For $K_1, K_2 \in \mathcal{K}^n$ and $\lambda_1, \lambda_2 \geq
0$, the support function of the Minkowski linear combination
$\lambda_1K_1+\lambda_2K_2$ is
\[h(\lambda_1K_1+\lambda_2K_2,\cdot)=\lambda_1h(K_1,\cdot)+\lambda_2h(K_2,\cdot).\]
By a theorem of Minkowski, the volume of a Minkowski linear
combination $\lambda_1K_1 + \ldots + \lambda_m K_m$ of convex
bodies $K_1, \ldots, K_m$ can be expressed as a homogeneous
polynomial of degree $n$:
\begin{equation} \label{mixed}
V(\lambda_1K_1 + \ldots +\lambda_m K_m)=\sum \limits_{i_1,\ldots,
i_n} V(K_{i_1},\ldots,K_{i_n})\lambda_{i_1}\cdots\lambda_{i_n}.
\end{equation}
The coefficients $V(K_{i_1},\ldots,K_{i_n})$ are called {\it mixed
volumes} of $K_{i_1}, \ldots, K_{i_n}$. These functionals are
nonnegative, symmetric and translation invariant. Clearly, their
diagonal form is ordinary volume, i.e., $V(K,\ldots,K)=V(K)$.

For $K, L \in \mathcal{K}^n$ and $0 \leq i \leq n - 1$, we write
$W_i(K,L)$ to denote the mixed volume
$V(K,\ldots,K,B,\ldots,B,L)$, where $K$ appears $n-i-1$ times and
the Euclidean unit ball $B$ appears $i$ times. The mixed volume
$W_i(K,K)$ will be written as $W_i(K)$ and is called the
\emph{$i$th quermassintegral} of $K$. The \emph{$i$th intrinsic
volume} $V_i(K)$ of $K$ is defined by
\begin{equation} \label{viwi}
\kappa_{n-i}V_i(K)={n \choose i}W_{n-i}(K),
\end{equation}
where $\kappa_n$ is the $n$-dimensional volume of the Euclidean
unit ball in $\mathbb{R}^n$. A special case of (\ref{mixed}) is
the classical Steiner formula for the volume of the outer
parallel body of $K$ at distance $\varepsilon > 0$:
\[V(K + \varepsilon B)=\sum \limits_{i=0}^n \varepsilon^i{n \choose i}W_i(K) = \sum \limits_{i=0}^n \varepsilon^{n-i}\kappa_{n-i}V_i(K).  \]

For $1 \leq i \leq n - 1$, the quermassintegral $W_{n-i}(K)$ of $K
\in \mathcal{K}^n$ can also be defined by
\begin{equation} \label{wigrin}
W_{n-i}(K)=\frac{\kappa_n}{\kappa_i}
\int_{\mathrm{Gr}_{i,n}}\mathrm{vol}_i(K|E)\,dE,
\end{equation}
where the integration is with respect to the rotation invariant
probability measure on $\mathrm{Gr}_{i,n}$.

\pagebreak

Let $\mathcal{K}^n_{\mathrm{o}}$ denote the set of convex bodies
in $\mathbb{R}^n$ with non-empty interior. \linebreak One of the
fundamental inequalities for mixed volumes is the general
Minkowski inequality: If $K, L \in \mathcal{K}^n_{\mathrm{o}}$
and $0 \leq i \leq n - 2$, then
\begin{equation} \label{genmink}
W_i(K,L)^{n-i} \geq W_i(K)^{n-i-1}W_i(L),
\end{equation}
with equality if and only if $K$ and $L$ are homothetic.

A consequence of the Minkowski inequality (\ref{genmink}) is the
Brunn--Minkowski inequality for quermassintegrals: If $K, L \in
\mathcal{K}^n_{\mathrm{o}}$ and $0 \leq i \leq n - 2$, then
\begin{equation} \label{quermassbm}
W_i(K+L)^{1/(n-i)} \geq W_i(K)^{1/(n-i)}+W_i(L)^{1/(n-i)},
\end{equation}
with equality if and only if $K$ and $L$ are homothetic.

For $K, K_1, \ldots, K_i \in \mathcal{K}^n$ and
$\mathbf{C}=(K_1,\ldots,K_i)$, let $V_i(K,\mathbf{C})$ denote the
mixed volume $V(K,\ldots,K,K_1,\ldots,K_i)$ with $n-i$ copies of
$K$. A further generalization of inequality (\ref{quermassbm})
(but without equality conditions) is the following: If $0 \leq i
\leq n-2$, $K, L, K_1, \ldots, K_i \in \mathcal{K}^n$ and
$\mathbf{C}=(K_1,...,K_i)$, then
\begin{equation} \label{mostgenbm}
V_i(K+L,\mathbf{C})^{1/(n-i)} \geq
V_i(K,\mathbf{C})^{1/(n-i)}+V_i(L,\mathbf{C})^{1/(n-i)}.
\end{equation}

A convex body $K \in \mathcal{K}^n_{\mathrm{o}}$ is also
determined up to translation by its surface area measure
$S_{n-1}(K,\cdot)$. Recall that for a Borel set $\omega \subseteq
S^{n-1}$, $S_{n-1}(K,\omega)$ is the $(n-1)$-dimensional
Hausdorff measure of the set of all boundary points of $K$ at
which there exists a normal vector of $K$ belonging to $\omega$.
The relation $S_{n-1}(\lambda
K,\cdot)=\lambda^{n-1}S_{n-1}(K,\cdot)$ holds for all $K \in
\mathcal{K}^n$ and every $\lambda \geq 0$. For $\vartheta \in
\mathrm{O}(n)$, we have $S_{n-1}(\vartheta K,\cdot)=\vartheta
S_{n-1}(K,\cdot)$, where $\vartheta S_{n-1}(K,\cdot)$ is the
image measure of $S_{n-1}(K,\cdot)$ under $\vartheta \in
\mathrm{O}(n)$.

The surface area measure $S_{n-1}(K,\cdot)$ of $K \in
\mathcal{K}^n$ satisfies the Steiner-type formula
\begin{equation} \label{steiner}
S_{n-1}(K + \varepsilon B,\cdot)=\sum \limits_{i=0}^{n-1}
\varepsilon^{n-1-i} {n - 1 \choose i}S_i(K,\cdot).
\end{equation}
The measure $S_i(K,\cdot)$ is called the {\it area measure of
order $i$} of $K \in \mathcal{K}^n$.

We conclude this section with an integral representation
connecting area measures and quermassintegrals: If $K, L \in
\mathcal{K}^n$ and $0 \leq i \leq n - 1$, then
\begin{equation} \label{wisi}
W_{n-1-i}(K,L)=\frac{1}{n} \int_{S^{n-1}}h(L,u)\,dS_i(K,u).
\end{equation}

\pagebreak

\centerline{\large{\bf{ \setcounter{abschnitt}{3}
\arabic{abschnitt}. Translation invariant valuations}}}

\vspace{0.7cm} \reseteqn \alpheqn \setcounter{koro}{0}

In this section, we collect some results from the theory of
real-valued translation invariant valuations. In particular, we
recall the definition of smooth valuations and the notion of
Crofton measures.

A function $\phi: \mathcal{K}^n \rightarrow \mathcal{A}$ into an
abelian semigroup $(\mathcal{A},+)$ is called a {\it valuation} if
\[\phi(K \cup L) + \phi(K \cap L)=\phi(K)+\phi(L),  \]
whenever $K, L, K \cup L \in \mathcal{K}^n$. The notion of
valuation as defined here is a classical concept from convex
geometry. However, we remark that Alesker has recently introduced
a broader notion of valuation in the more general setting of
smooth manifolds, see \textbf{\cite{Alesker06b, Alesker06a,
Alesker07, AlFu, Bernig08b}}.

A valuation $\phi$ is called {\it translation invariant} if
$\phi(K+x)=\phi(K)$ for all $x \in \mathbb{R}^n$ and $K \in
\mathcal{K}^n$. We denote the vector space of continuous
translation invariant {\it real-valued} valuations by
$\mathbf{Val}$ and we write $\mathbf{Val}_i$ for its subspace of
all valuations of degree $i$. A valuation $\phi \in \mathbf{Val}$
is said to be even (resp.\ odd) if
$\phi(-K)=(-1)^{\alpha}\phi(K)$ with $\alpha = 0$ (resp.\ $\alpha
= 1$) for all $K \in \mathcal{K}^n$. We write $\mathbf{Val}_i^+
\subseteq \mathbf{Val}_i$ for the subspace of even valuations of
degree $i$ and $\mathbf{Val}_i^-$ to denote the space of odd
valuations of degree $i$, respectively.

The following result was obtained by McMullen:

\begin{satz} \label{mcmullen} \emph{(McMullen \textbf{\cite{McMullen77}})}
\[\mathbf{Val} = \bigoplus \limits_{i=0}^n \mathbf{Val}_i^+ \oplus \mathbf{Val}_i^-.   \]
\end{satz}

It follows from Theorem \ref{mcmullen} that the space
$\mathbf{Val}$ becomes a Banach space under the norm
\[\|\phi \| = \sup \{|\phi(K)|: K \subseteq B \}.  \]

\vspace{0.3cm}

\noindent {\bf Examples:}
\begin{enumerate}
\item[(a)] It is easy to see that the space $\mathbf{Val}_0$ is one-dimensional
and is spanned by the Euler characteristic $\chi$ (recall that
$\chi(K)=1$ for all $K \in \mathcal{K}^n$).
\item[(b)] Hadwiger \textbf{\cite[\textnormal{p.\
79}]{hadwiger51}} has shown that $\mathbf{Val}_n$ is also
one-dimensional and is spanned by ordinary volume $V$.
\item[(c)] For $i \in \{0, \ldots, n\}$, let us fix convex bodies $\mathbf{C}=(L_1, \ldots, L_{i})$. The mixed volume
$V_i(K,\mathbf{C})$ belongs to $\mathbf{Val}_{n-i}$.
\end{enumerate}

\pagebreak

The group $\mathrm{GL}(n)$ has a natural continuous
representation $\rho$ on the \linebreak Banach space
$\mathbf{Val}$: For every $A \in \mathrm{GL}(n)$ and $K \in
\mathcal{K}^n$,
\[(\rho(A)\phi)(K)=\phi(A^{-1}K), \qquad \phi \in \mathbf{Val}. \]
Note that the subspaces $\mathbf{Val}_i^{\pm} \subseteq
\mathbf{Val}$ are invariant under this $\mathrm{GL}(n)$ action.
The Irreducibility Theorem of Alesker states the following:

\begin{satz} \label{irr} \emph{(Alesker \textbf{\cite{Alesker01}})} The
natural representation of $\mathrm{GL}(n)$ on
$\mathbf{Val}_i^{\pm}$ is irreducible for any $i \in \{0, \ldots,
n\}$.
\end{satz}

The Irreducibility Theorem directly implies a conjecture by
McMullen that the linear combinations of mixed volumes are dense
in $\mathbf{Val}$ (see \textbf{\cite{Alesker01}}).

In the following we will further illustrate the strength of the
Irreducibility Theorem by constructing an alternative description
of translation invariant even valuations: Assume that $1 \leq i
\leq n - 1$. For any finite Borel measure $\mu$ on
$\mathrm{Gr}_{i,n}$ define an even valuation $\mathrm{A}_i\mu \in
\mathbf{Val}_i^+$ by
\begin{equation} \label{crofton}
(\mathrm{A}_i\mu)(K)=\int_{\mathrm{Gr}_{i,n}}\mathrm{vol}_i(K|E)\,d\mu(E).
\end{equation}
Clearly, the image of the map $\mathrm{A}_i$ is a
$\mathrm{GL}(n)$ invariant subspace of $\mathbf{Val}_i^+$.
Therefore, by Theorem \ref{irr}, this image is dense in
$\mathbf{Val}_i^+$.

\vspace{0.3cm}

\noindent {\bf Definition} \emph{A finite Borel measure $\mu$ on
$\mathrm{Gr}_{i,n}$, $1 \leq i \leq n - 1$, is called a Crofton
measure for the valuation $\phi \in \mathbf{Val}_i^+$ if
$\mathrm{A}_i\mu = \phi$.}

\vspace{0.3cm}

The classical Crofton formula is a result from the early days of
integral geometry relating the length of a curve in the plane to
the expected number of intersection points with random lines.
Higher-dimensional generalizations have become known as linear
kinematic formulas. For more information on Crofton formulas the
reader is referred to the recent book by Schneider and Weil
\textbf{\cite{schnweil}}. Further details on Crofton measures of
valuations (not necessarily translation invariant) can be found in
\textbf{\cite{Bernig07}}.

In the following it will be important for us to work with a
subset of valuations in $\mathbf{Val}_i^+$ which admit a Crofton
formula (\ref{crofton}).

\vspace{0.3cm}

\noindent {\bf Definition} \emph{A valuation $\phi \in
\mathbf{Val}$ is called smooth if the map $\mathrm{GL}(n)
\rightarrow \mathbf{Val}$ defined by $A \mapsto \rho(A)\phi$ is
infinitely differentiable.}

\vspace{0.3cm}

The notion of smooth valuations is a special case of the more
general concept of smooth vectors in a representation space (see,
e.g., \textbf{\cite[\textnormal{p.\ 31}]{wallach1}}).

\pagebreak

We denote the space of smooth translation invariant valuations by
$\mathbf{Val}^{\infty}$ and we write $\mathbf{Val}_i^{\pm,\infty}$
for the subspace of smooth valuations in $\mathbf{Val}_i^{\pm}$.
From representation theory it is well known (cf.\
\textbf{\cite[\textnormal{p.\ 32}]{wallach1}}) that the set of
smooth valuations $\mathbf{Val}_i^{\pm,\infty}$ is a dense
$\mathrm{GL}(n)$ invariant subspace of $\mathbf{Val}_i^{\pm}$ and
one easily deduces the following decomposition:
\begin{equation} \label{decsmooth}
\mathbf{Val}^{\infty} = \bigoplus \limits_{i=0}^n
\mathbf{Val}_i^{+,\infty} \oplus \mathbf{Val}_i^{-,\infty}.
\end{equation}

Now consider the restriction of the map $\mathrm{A}_i$, $1 \leq i
\leq n - 1$, defined in (\ref{crofton}) to smooth functions:
\begin{equation*}
(\mathrm{A}_if)(K)=\int_{\mathrm{Gr}_{i,n}}\mathrm{vol}_i(K|E)f(E)\,dE,
\qquad f \in C^{\infty}(\mathrm{Gr}_{i,n}).
\end{equation*}
Clearly, the valuation $\mathrm{A}_if$ is smooth, i.e.,
$\mathrm{A}_if \in \mathbf{Val}_i^{+,\infty}$. Moreover, it
follows from a deep result of Alesker and Bernstein
\textbf{\cite{AlBern}} that any smooth translation invariant and
even valuation admits such a Crofton formula. In order to explain
this fact we need the cosine transform on Grassmannians.

Assume that $1 \leq i \leq n - 1$. For two subspaces $E, F \in
\mathrm{Gr}_{i,n}$, the cosine of the angle between $E$ and $F$
is defined by
\[|\cos(E,F)|=\mathrm{vol}_i(\mathrm{Pr}_F(M)),  \]
where $M$ is any subset of $E$ with $\mathrm{vol}_i(M)=1$ and
$\mathrm{Pr}_F$ denotes the orthogonal projection onto $F$. (This
definition does not depend on the choice of $M \subseteq E$.) The
\emph{cosine transform} $\mathrm{C}_i: C(\mathrm{Gr}_{i,n})
\rightarrow C(\mathrm{Gr}_{i,n})$ is defined by
\[(\mathrm{C}_if)(F)=\int_{\mathrm{Gr}_{i,n}} |\cos(E,F)|f(E)\,dE.  \]

Alesker and Bernstein established a fundamental connection
between the range of the cosine transform and even translation
invariant valuations. This result is based on an imbedding
$\mathrm{K}_i: \mathbf{Val}_i^+ \rightarrow C(\mathrm{Gr}_{i,n})$
due to Klain: For $\phi \in \mathbf{Val}_i^+$ and every $E \in
\mathrm{Gr}_{i,n}$, consider the restriction $\phi_E$ of $\phi$ to
convex bodies in $E$. This is a continuous translation invariant
valuation of degree $i$ in $E$. Hence, by a result of Hadwiger
\textbf{\cite[\textnormal{p.\ 79}]{hadwiger51}}, $\phi_E =
g(E)\,\mathrm{vol}_i$, where $g(E)$ is a constant depending on
$E$. The map $\mathrm{K}_i: \mathbf{Val}_i^+ \rightarrow
C(\mathrm{Gr}_{i,n})$, defined by $\mathrm{K}_i\phi=g$, turns out
to be injective by a result of Klain \textbf{\cite{klain00}}. The
function $g$ is called the \emph{Klain function} of the valuation
$\phi$.

\pagebreak

The Alesker--Bernstein theorem can be stated as follows:

\begin{satz} \label{albern} \emph{(Alesker and Bernstein \textbf{\cite{AlBern}})}  Suppose that $1 \leq i \leq n - 1$.
The image of the Klain imbedding $\mathrm{K}_i:
\mathbf{Val}_i^{+,\infty} \rightarrow
C^{\infty}(\mathrm{Gr}_{i,n})$ coincides with the image of the
cosine transform $\mathrm{C}_i: C^{\infty}(\mathrm{Gr}_{i,n})
\rightarrow C^{\infty}(\mathrm{Gr}_{i,n})$.
\end{satz}

We remark that this version of the Alesker--Bernstein theorem is
obtained from the main results in \textbf{\cite{AlBern}} by an
application of the Casselman--Wallach theorem
\textbf{\cite{cassel}} (cf.\ \textbf{\cite[\textnormal{p.\
72}]{Alesker03}}).

For $1 \leq i \leq n - 1$, let $\mathbf{T}_i^{\infty}$ denote the
image of smooth functions on $\mathrm{Gr}_{i,n}$ under the cosine
transform $\mathrm{C}_i$. It is well known that $\mathrm{C}_i$ is
not injective for $2 \leq i \leq n - 2$. However, since
$\mathrm{C}_i$ is selfadjoint its restriction to
$\mathbf{T}_i^{\infty}$ has trivial kernel. Moreover, from an
application of the Casselman--Wallach theorem
\textbf{\cite{cassel}} to the main result of
\textbf{\cite{AlBern}}, Alesker \textbf{\cite[\textnormal{p.\
73}]{Alesker03}} deduced that
\begin{equation} \label{iso}
\mathrm{C}_i(\mathbf{T}_i^{\infty})=\mathbf{T}_i^{\infty}.
\end{equation}

Now suppose that $F \in \mathrm{Gr}_{i,n}$. Then, for any $f \in
C(\mathrm{Gr}_{i,n})$ and any convex body $K \subseteq F$,
\begin{equation} \label{ciki}
(\mathrm{A}_if)(K)=\mathrm{vol}_i(K)\int_{\mathrm{Gr}_{i,n}}
|\cos(E,F)|f(E)\,dE.
\end{equation}
Consequently, the Klain function of the valuation $\mathrm{A}_if$
is equal to the cosine transform $\mathrm{C}_if$ of $f$. Thus, we
obtain from Theorem \ref{albern} and (\ref{iso}):

\begin{koro} \label{crofkoro} For any valuation $\phi \in
\mathbf{Val}_i^{+,\infty}$, there exists a unique smooth measure
$\mu \in \mathbf{T}_i^{\infty}$ such that $\mu$ is a Crofton
measure for $\phi$.

\end{koro}

We conclude this section with a commutative diagram for the
crucial isomorphisms needed in the following:

\[
\begin{xy}
 \xymatrix{
       &     &  \mathbf{Val}_i^{+,\infty} \ar[dd]^{\textnormal{\normalsize $\mathrm{K}_i$}}  \\
       &     &     \\
      \mathbf{T}_i^{\infty}  \ar[uurr]^{\textnormal{\normalsize $\mathrm{A}_i$}}    \ar[rr]^{\textnormal{\normalsize $\mathrm{C}_i$}} &     & \,\, \mathbf{T}_i^{\infty}
  }
\end{xy}
\]

\pagebreak

\centerline{\large{\bf{ \setcounter{abschnitt}{4}
\arabic{abschnitt}. Convolutions}}}

\vspace{0.7cm} \reseteqn \alpheqn \setcounter{koro}{0}

Here we recall the basic notion of convolution on the compact Lie
group $\mathrm{O}(n)$ and the homogeneous spaces
$\mathrm{O}(n)/\mathrm{O}(n-1)$ and $\mathrm{O}(n)/\mathrm{O}(i)
\times \mathrm{O}(n - i)$. At the end of this section, we
establish an auxiliary result which is critical in the proof of
Theorem \ref{thm3}. As a general reference for this section we
recommend the article by Grinberg and Zhang
\textbf{\cite{grinbergzhang99}}.

Let $C(\mathrm{O}(n))$ denote the space of continuous functions on
$\mathrm{O}(n)$ with the uniform topology. In this article all
measures are signed finite Borel measures. For $f \in
C(\mathrm{O}(n))$ and a measure $\mu$ on $\mathrm{O}(n)$, the
canonical pairing is
\begin{equation*}
\langle \mu,f\rangle=\langle f,\mu
\rangle=\int_{\mathrm{O}(n)}f(\vartheta)\,d\mu(\vartheta).
\end{equation*}
We will frequently identify a continuous function $f \in
C(\mathrm{O}(n))$ with the absolutely continuous measure (with
respect to Haar probability measure on $\mathrm{O}(n)$) with
density $f$. The canonical pairing is then consistent with the
usual inner product on $C(\mathrm{O}(n))$.

For $\vartheta \in \mathrm{O}(n)$, the left translation $\vartheta
f$ of $f \in C(\mathrm{O}(n))$ is defined by
\begin{equation*}
\vartheta f(\eta)=f(\vartheta^{-1}\eta).
\end{equation*}
For a measure $\mu$ on $\mathrm{O}(n)$, we set
\begin{equation*}
\langle \vartheta \mu,f\rangle=\langle
\mu,\vartheta^{-1}f\rangle, \qquad f \in C(\mathrm{O}(n)).
\end{equation*}
Then $\vartheta \mu$ is just the image measure of $\mu$ under the
rotation $\vartheta$.

For $f \in C(\mathrm{O}(n))$, the function $\widehat{f} \in
C(\mathrm{O}(n))$ is defined by
\begin{equation*}
\widehat{f}(\vartheta)=f(\vartheta^{-1}).
\end{equation*}
For a measure $\mu$ on $\mathrm{O}(n)$, we define the measure
$\widehat{\mu}$ by
\begin{equation*}
\langle \widehat{\mu},f\rangle=\langle \mu,\widehat{f}\rangle,
\qquad f \in C(\mathrm{O}(n)).
\end{equation*}

For $f, g \in C(\mathrm{O}(n))$, the convolution $f \ast g \in
C(\mathrm{O}(n))$ is defined by
\[(f \ast g)(\eta)=\int_{\mathrm{O}(n)} f(\eta \vartheta^{-1})g(\vartheta)\,d\vartheta
=\int_{\mathrm{O}(n)}f(\vartheta)g(\vartheta^{-1}\eta)\,d\vartheta,\]
where integration is with respect to the Haar probability measure
on $\mathrm{O}(n)$.

\noindent For a measure $\mu$ on $\mathrm{O}(n)$ and a function
$f \in C(\mathrm{O}(n))$, the convolutions $\mu \ast f \in
C(\mathrm{O}(n))$ and $f \ast \mu \in C(\mathrm{O}(n))$ are
defined by
\begin{equation} \label{convmeas}
(f \ast \mu)(\eta)=\int_{\mathrm{O}(n)} f(\eta
\vartheta^{-1})\,d\mu(\vartheta), \quad (\mu \ast
f)(\eta)=\int_{\mathrm{O}(n)}\vartheta f(\eta)\,d\mu(\vartheta).
\end{equation}
From this definition, it follows that $f \ast \mu$ and $\mu \ast
f$ are $C^{\infty}$ if $f \in C^{\infty}(\mathrm{O}(n))$.

We emphasize that, if $\mu$ is a measure on $\mathrm{O}(n)$,
then, by (\ref{convmeas}), for all $f \in C(\mathrm{O}(n))$ and
every $\vartheta \in \mathrm{O}(n)$,
\begin{equation} \label{roteq}
(\vartheta f) \ast \mu = \vartheta(f \ast \mu).
\end{equation}
Thus, the convolution from the right gives rise to operators on
$C(\mathrm{O}(n))$ which intertwine orthogonal transformations.

Using (\ref{convmeas}), it is also easy to verify that for $f, g
\in C(\mathrm{O}(n))$ and a measure $\sigma$ on $\mathrm{O}(n)$,
\begin{equation} \label{adjoint}
\langle g \ast \sigma ,f \rangle = \langle g,f \ast
\widehat{\sigma}\rangle.
\end{equation}
This identity allows us to define the convolution $\mu \ast
\sigma$ of two measures $\mu, \sigma$ on $\mathrm{O}(n)$ by
\begin{equation*}
\langle\mu \ast \sigma,f \rangle=\langle \sigma,\widehat{\mu} \ast
f\rangle = \langle \mu,f\ast \widehat{\sigma}\rangle, \qquad f
\in C(\mathrm{O}(n)).
\end{equation*}

It is easy to check that the convolution of functions and measures
on $\mathrm{O}(n)$ defined in this way is associative. However,
the convolution is in general not commutative. If $\mu, \sigma$
are measures on $\mathrm{O}(n)$, then
\begin{equation} \label{commute}
\widehat{\mu \ast \sigma}=\widehat{\sigma} \ast \widehat{\mu}.
\end{equation}

For the following Lemma see \textbf{\cite[\textnormal{p.\
85}]{grinbergzhang99}}.

\begin{lem} \label{approxlem} Let $\mu, \mu_m$, $m \in \mathbb{N}$, be measures on
$\mathrm{O}(n)$ and let $f \in C(\mathrm{O}(n))$. If $\mu_m
\rightarrow \mu$ weakly, then $f \ast \mu_m \rightarrow f \ast
\mu$ and $\mu_m \ast f \rightarrow \mu \ast f$ uniformly.
\end{lem}

In the following we extend the definition of convolutions to
functions and measures on the homogeneous spaces
\[S^{n-1}=\mathrm{O}(n)/\mathrm{O}(n-1) \qquad \mbox{and} \qquad \mathrm{Gr}_{i,n}=\mathrm{O}(n)/\mathrm{O}(i) \times \mathrm{O}(n-i).  \]
In order to treat both cases simultaneously let $H$ denote a
closed subgroup of $\mathrm{O}(n)$. We consider the compact
homogeneous space $\mathrm{O}(n)/H$.

Let $\pi: \mathrm{O}(n) \rightarrow \mathrm{O}(n)/H$ be the
canonical projection and write $\pi(\vartheta)= \bar{\vartheta}$.
If $e \in \mathrm{O}(n)$ denotes the identity map, then $H$ is the
stabilizer in $\mathrm{O}(n)$ of $\bar{e} \in \mathrm{O}(n)/H$
and we have $\bar{\vartheta}=\vartheta \bar{e}$ for every
$\vartheta \in \mathrm{O}(n)$.

\pagebreak

Every continuous function $f$ on $\mathrm{O}(n)/H$ gives rise to
a continuous right $H$-invariant function $\breve{f}=f \circ \pi$
on $\mathrm{O}(n)$. Conversely, every $f \in C(\mathrm{O}(n))$
induces a continuous function $\bar{f}$ on $\mathrm{O}(n)/H$,
defined by
\[\bar{f}(\bar{\eta})=\int_{H}f(\eta\vartheta)\,d\vartheta. \]
If $f \in C(\mathrm{O}(n))$ is right $H$-invariant, then $f =
\bar{f} \circ \pi$. Therefore, the subspace of right $H$-invariant
functions in $C(\mathrm{O}(n))$ is isomorphic to
$C(\mathrm{O}(n)/H)$.

For a measure $\mu$ on $\mathrm{O}(n)/H$, we define the measure
$\breve{\mu}$ on $\mathrm{O}(n)$ by
\[\langle\, \breve{\mu}, f \rangle = \langle\, \mu , \bar{f} \rangle.\]
In this way we also obtain a one-to-one correspondence between
measures on $\mathrm{O}(n)/H$ and right $H$-invariant measures on
$\mathrm{O}(n)$.

The convolutions of functions and measures on $\mathrm{O}(n)/H$
can be defined via their identification with right $H$-invariant
functions and measures on $\mathrm{O}(n)$. For example, the
convolution $f \ast \mu \in C(\mathrm{O}(n)/H)$ of $f \in
C(\mathrm{O}(n)/H)$ with a measure $\mu$ on $\mathrm{O}(n)/H$ is
defined by
\begin{equation} \label{convhom}
(f \ast \mu)(\bar{\eta})=(\breve{f} \ast \breve{\mu})(\eta)
=\int_{\mathrm{O}(n)}f(\eta\vartheta^{-1}\bar{e})\,d\breve{\mu}(\vartheta).
\end{equation}
In the same way we can define convolutions between different
homogeneous spaces: Let $H_1, H_2$ be two closed subgroups of
$\mathrm{O}(n)$. If $f \in C(\mathrm{O}(n)/H_1)$ and $g \in
C(\mathrm{O}(n)/H_2)$, then $\breve{f} \ast \breve{g}$ defines a
continuous right $H_2$-invariant function on $\mathrm{O}(n)$ and
thus can be identified with a continuous function on
$\mathrm{O}(n)/H_2$.

It follows from (\ref{convhom}) that the Dirac measure
$\delta_{\bar{e}}$ on $\mathrm{O}(n)/H$ is the unique
rightneutral element for the convolution on $\mathrm{O}(n)/H$. If
$f \in C(\mathrm{O}(n))$, then
\begin{equation} \label{dirac1}
f \ast \delta_{\bar{e}} = \bar{f}
\end{equation}
is right $H$-invariant and
\begin{equation} \label{dirac2}
\delta_{\bar{e}} \ast f  =\int_{H}\vartheta f\,d\vartheta
\end{equation}
defines a left $H$-invariant function on $\mathrm{O}(n)$.

An essential role among functions (and measures) on
$\mathrm{O}(n)$ play the \linebreak $H$-biinvariant functions.
They can be identified with functions on $\mathrm{O}(n)/H$ with
the property that $\vartheta f = f$ for every $\vartheta \in H$.
We call a $H$-invariant function on $\mathrm{O}(n)/H$ {\it
zonal}. If $f, g \in C(\mathrm{O}(n)/H)$, then, by (\ref{dirac1})
and (\ref{dirac2}),
\begin{equation} \label{zonsuff}
f \ast g = (f \ast \delta_{\bar{e}}) \ast g = f \ast
(\delta_{\bar{e}} \ast g).
\end{equation}

\pagebreak

\noindent Consequently, for convolutions from the right on
$\mathrm{O}(n)/H$, it is sufficient to consider zonal functions
and measures.

If $f \in C(\mathrm{O}(n))$ is $H$-biinvariant (or, equivalently,
$f \in C(\mathrm{O}(n)/H)$ is zonal), then the function
$\widehat{f}$ is also $H$-biinvariant and thus can be identified
with a zonal function on $\mathrm{O}(n)/H$.

It is trivial to verify that if $f \in C(S^{n-1})$ is zonal, then
$\widehat{f}=f.$ The corresponding result for the Grassmannian
will be crucial in the proof of Theorem \ref{thm3}.

\begin{lem} \label{selfad} If $f \in C(\mathrm{Gr}_{i,n})$ is $\mathrm{O}(i) \times \mathrm{O}(n-i)$ invariant, then
\[\widehat{f} = f.  \]
\end{lem}
{\it Proof.} Let $H$ denote the subgroup $\mathrm{O}(i) \times
\mathrm{O}(n - i)$. We identify $f$ with a $H$-biinvariant
function on $\mathrm{O}(n)$. We will show that for any $\vartheta
\in \mathrm{O}(n)$,
\begin{equation} \label{fdach}
H\,\vartheta^{-1}H=H\,\vartheta\,H.
\end{equation}
Clearly, (\ref{fdach}) implies that $\widehat{f}=f$.

For the following short proof of (\ref{fdach}) the author is
obliged to S. Alesker. Replacing $i$ with $n - i$ if necessary,
we may assume that $i \leq n/2$. Fix an orthonormal basis $\{b_1,
\ldots, b_{n}\}$ of $\mathbb{R}^n$ and let $I_n$ denote the
identity matrix. Consider the torus $T \subseteq \mathrm{O}(n)$
consisting of rotations of the form
\[\left ( \begin{array}{l|c}
 \begin{array}{ccc|ccc} \cos \alpha_1 &  &  & - \sin \alpha_1 & & \\  & \ddots & & & \ddots
 &\\ & &  \cos \alpha_i &  &  & - \sin \alpha_i \\ \hline  \sin \alpha_1 &  &  & \cos \alpha_1 & & \\  & \ddots & & & \ddots
 &\\ & &  \sin \alpha_i &  &  & \cos \alpha_i
\end{array} &  \\ \hline & I_{n-2i}
\end{array}
 \right ).   \]
It is well known (and can be generalized appropriately to the
setting of Riemannian symmetric pairs (see, e.g.,
\textbf{\cite[\textnormal{Chapter II}]{takeuchi}})) that
\[\mathrm{O}(n)=H\,T\,H.   \]
Consequently, we can assume for the proof of (\ref{fdach}) that
$\vartheta \in T$. But now it is straightforward to verify that
$\vartheta^{-1}=J\,\vartheta\,J$, where $J \in H$ is given by
\[J=\left (\begin{array}{cc} -I_i & 0 \\ 0 & I_{n-i}   \end{array}   \right ).  \]
\hfill $\blacksquare$

\pagebreak

As a consequence of (\ref{commute}) and Lemma \ref{selfad}, we
note that the convolution of zonal functions (and measures) on
$S^{n-1}$ and $\mathrm{Gr}_{i,n}$ is abelian.

Another important ingredient in the proof of Theorem \ref{thm3}
are spherical approximate identities. Let $B_m(\bar{e})$ be the
open geodesic ball of radius $\frac{1}{m}$ at $\bar{e} \in
S^{n-1}$, where $m \in \mathbb{N}$ is sufficiently large. A
sequence $f_m$, $m \in \mathbb{N}$, of nonnegative $C^{\infty}$
functions on $S^{n-1}$ is called a {\it spherical approximate
identity} if, for each $m$, the following two conditions are
satisfied:
\begin{enumerate}
\item[(i)] $\int_{S^{n-1}}f_m(u)\,du=1$;
\item[(ii)] $f_m(u)=0$ if $u \not \in B_m(\bar{e})$.

\end{enumerate}

The existence of spherical approximate identities follows from
standard techniques similar to the construction of partitions of
unity on manifolds (cf.\ \textbf{\cite[\textnormal{p.\
84}]{grinbergzhang99}}). We conclude this section with a
well-known auxiliary result, see, e.g.,
\textbf{\cite[\textnormal{Lemma 2.5}]{grinbergzhang99}}.

\begin{lem} \label{approxid} If $f_m$, $m \in \mathbb{N}$, is a spherical approximate identity, then
\begin{enumerate}
\item[(a)] $\lim_{m \rightarrow \infty} g \ast f_m = g$ uniformly
for every $g \in C(S^{n-1})$;
\item[(b)] $\lim_{m \rightarrow \infty} \mu \ast f_m = \mu$
weakly for every measure $\mu$ on $S^{n-1}$.
\end{enumerate}
\end{lem}

\vspace{1cm}

\centerline{\large{\bf{ \setcounter{abschnitt}{5}
\arabic{abschnitt}. Minkowski valuations}}}

\vspace{0.7cm} \reseteqn \alpheqn \setcounter{koro}{0}

In this section we collect the background material on translation
invariant Minkowski valuations. We also extend the definition of
smooth valuations to translation invariant Minkowski valuations
which are $\mathrm{O}(n)$ equivariant.

A Minkowski valuation $\Phi: \mathcal{K}^n \rightarrow
\mathcal{K}^n$ is called $\mathrm{O}(n)$ {\it equivariant} if for
all $K \in \mathcal{K}^n$ and every $\vartheta \in \mathrm{O}(n)$,
\[\Phi(\vartheta K)= \vartheta \Phi K.  \]

We denote by $\mathbf{MVal}$ the set of continuous translation
invariant Minkowski valuations which are $\mathrm{O}(n)$
equivariant and we write $\mathbf{MVal}_i^{(+)}$ for its subset of
all (even) Minkowski valuations of degree $i$. (This slight abuse
of notation should not lead to confusion because in the following
all Minkowski valuations will be translation invariant {\it and}
$\mathrm{O}(n)$ equivariant.)

Since Minkowski valuations arise naturally, like the projection
and the \linebreak difference operator, from data about
projections and sections of convex bodies, they form an integral
part of geometric tomography. In the following we give a few
well-known examples of Minkowski valuations in $\mathbf{MVal}$
(for additional examples, see, e.g., \textbf{\cite{kiderlen05,
Ludwig:Minkowski, Schu08}}):

\vspace{0.3cm}

\noindent {\bf Examples:}
\begin{enumerate}
\item[(a)] For $i \in \{0, \ldots, n\}$, define $\Lambda_i:
\mathcal{K}^n \rightarrow \mathcal{K}^n$ by
$\Lambda_i(K)=V_i(K)B$. Clearly, we have $\Lambda_i \in
\mathbf{MVal}_i^+$.
\item[(b)] For $i \in \{1, \ldots, n - 1\}$, the $i$th projection
operator $\Pi_i: \mathcal{K}^n \rightarrow \mathcal{K}^n$ is an
element of $\mathbf{MVal}_i^+$.
\item[(c)] For $i \in \{1, \ldots, n - 1\}$, the (normalized) $i$th
mean section operator $\mathrm{M}_i \in \mathbf{MVal}_{n+1-i}$,
introduced by Goodey and Weil \textbf{\cite{goodeyweil92}}, is
defined by
\[h(\mathrm{M}_iK,\cdot) = \int_{\mathrm{AGr}_{i,n}} h(K \cap E, \cdot)\,d\mu_i(E) - h(\{z_{n+1-i}\},\cdot).  \]
Here, $\mathrm{AGr}_{i,n}$ is the affine Grassmannian of
$i$-dimensional planes in $\mathbb{R}^n$, $\mu_i$ is its
(suitably normalized) motion invariant measure and $z_i(K)$
denotes the $i$th moment vector of $K$ (see
\textbf{\cite[\textnormal{p.\ 304}]{schneider93}}).
\end{enumerate}

\vspace{0.2cm}

From Theorem \ref{mcmullen}, one can deduce the following
decomposition result for Minkowski valuations (cf.\
\textbf{\cite[\textnormal{p.\ 12}]{schnschu}}):

\begin{lem} \label{lemdec} Suppose that $\Phi \in \mathbf{MVal}^{(+)}$. Then there
are constants $c_0, c_n \geq 0$ such that for every $K \in
\mathcal{K}^n$,
\[h(\Phi K,\cdot) = c_0 + \sum \limits_{i=1}^{n-1} g_i(K,\cdot) + c_nV(K),  \]
where the (even) function $g_i(K,\cdot) \in C(S^{n-1})$ has the
following properties:
\begin{enumerate}
\item[(a)] The map $K \mapsto g_i(K,\cdot)$ is a continuous
translation invariant (even) valuation of degree $i$.
\item[(b)] For every $\vartheta \in \mathrm{O}(n)$ and $K \in \mathcal{K}^n$, we have
$g_i(\vartheta K,u)=g_i(K,\vartheta^{-1}u)$.
\end{enumerate}
\end{lem}

It is not known, at this point, whether, for every $K \in
\mathcal{K}^n$, each function \linebreak $g_i(K,\cdot)$ is the
support function of a convex body. Hence, the following \linebreak
important problem is still open:

\pagebreak

\noindent {\bf Open Problem} \emph{Suppose that $\Phi \in
\mathbf{MVal}$. Is there a (unique) representation of $\Phi$ of
the form
\[\Phi = \Phi_0 + \ldots + \Phi_n,  \]
where $\Phi_i \in \mathbf{MVal}_i$?}

\vspace{0.3cm}

Suppose that $\Phi \in \mathbf{MVal}$. We define a {\it
real-valued} translation invariant valuation $\varphi \in
\mathbf{Val}$ by
\begin{equation} \label{real}
\varphi(K)=h(\Phi K,\bar{e}), \qquad K \in \mathcal{K}^n.
\end{equation}
Since $\Phi$ is $\mathrm{O}(n)$ equivariant, we have for
$\bar{\eta} \in S^{n-1}$,
\begin{equation} \label{heta}
h(\Phi K,\bar{\eta})= h(\Phi K,\eta \bar{e})= h(\Phi(\eta^{-1}
K), \bar{e})=\varphi(\eta^{-1}K).
\end{equation}
Consequently, the real-valued valuation $\varphi$ uniquely
determines the Minkowski valuation $\Phi$. We call the valuation
$\varphi \in \mathbf{Val}$ defined by (\ref{real}) the {\it
associated real-valued valuation} of $\Phi \in \mathbf{MVal}$.

We can now extend the notion of smooth real-valued valuations to
Minkowski valuations in $\mathbf{MVal}$.

\vspace{0.3cm}

\noindent {\bf Definition} \emph{A Minkowski valuation $\Phi \in
\mathbf{MVal}$ is called smooth if its \linebreak associated
real-valued valuation $\varphi \in \mathbf{Val}$ is smooth.}

\vspace{0.3cm}

We denote by $\mathbf{MVal}^{\infty}$ the subset of smooth
Minkowski valuations in $\mathbf{MVal}$ and we write
$\mathbf{MVal}_i^{(+),\infty}$ for the subset of smooth (even)
Minkowski valuations in $\mathbf{MVal}_i^{(+)}$.

\vspace{0.2cm}

A description of Minkowski valuations in $\mathbf{MVal}_1$ was
recently obtained by Kiderlen \textbf{\cite[\textnormal{Theorem
1.3}]{kiderlen05}} (extending previous results by Schneider
\textbf{\cite{schneider74}}). Here, we state a version of
Kiderlen's result for smooth Minkowski valuations:

\begin{satz} \label{kiderlen} \textnormal{(Kiderlen \textbf{\cite{kiderlen05}})} Suppose that $\Phi \in \mathbf{MVal}_1^{\infty}$.
Then there exists a unique smooth zonal measure $\mu$ on $S^{n-1}$
such that for every $K \in \mathcal{K}^n$,
\begin{equation} \label{kid}
h(\Phi K,\cdot) = h(K,\cdot) \ast \mu.
\end{equation}
Moreover, the Minkowski valuation $\Phi$ is even if and only if
$\mu$ is even.
\end{satz}

Since $h(K,u) + h(-K,u) = \mathrm{vol}_1(K|u)$, we remark that for
$\Phi \in \mathbf{MVal}_1^{+,\infty}$, representation (\ref{kid})
is equivalent to
\[h(\Phi K,\cdot) = \mathrm{vol}_1(K|\,\cdot\,) \ast \mu.   \]

Note that Theorem \ref{kiderlen} is not a complete
characterization of Minkowski valuations in
$\mathbf{MVal}_1^{\infty}$ but only a representation result. It is
not known which zonal measures on $S^{n-1}$ define a Minkowski
valuation by (\ref{kid}). However, the following conjecture
appears implicitly in \textbf{\cite{kiderlen05}}:

\vspace{0.3cm}

\noindent {\bf Conjecture} \emph{A map $\Phi: \mathcal{K}^n
\rightarrow \mathcal{K}^n$ is a Minkowski valuation in
$\mathbf{MVal}_1$ if and only if there exists a zonal measure
$\mu$ on $S^{n-1}$ which is non-negative up to addition of a
function $u \mapsto x \cdot u$, $x \in \mathbb{R}^n$, such that
for every $K \in \mathcal{K}^n$,}
\[h(\Phi K,\cdot) = h(K,\cdot) \ast \mu.\]

\vspace{0.2cm}

Recently, the author established a result corresponding to
Theorem \ref{kiderlen} for Minkowski valuations in
$\mathbf{MVal}_{n-1}$ (see \textbf{\cite[\textnormal{Theorem 1.2
\& 1.3}]{Schu06a}}):

\begin{satz} \label{schuster} Suppose that $\Phi \in \mathbf{MVal}_{n-1}$.
Then there exists a (unique) zonal function $g \in C(S^{n-1})$
such that for every $K \in \mathcal{K}^n$,
\begin{equation} \label{schu}
h(\Phi K,\cdot) = S_{n-1}(K,\cdot) \ast g.
\end{equation}
Moreover, the Minkowski valuation $\Phi$ is even if and only if $
g = h(L,\cdot)$, for some origin-symmetric body of revolution $L
\in \mathcal{K}^n$.
\end{satz}

It is well known that $\mathrm{vol}_{n-1}(K|\,\cdot\,)$ for $K
\in \mathcal{K}^n$ is (up to a constant factor) given by the
cosine transform of $S_{n-1}(K,\cdot)$. More precisely,
\[\mathrm{vol}_{n-1}(K|u^{\bot})=\frac{1}{2}\int_{S^{n-1}}|u\cdot v|\,dS_{n-1}(K,v)=(S_{n-1}(K,\cdot)\ast \mbox{$\frac{1}{2}$}|\bar{e}\cdot\,.\,|)(u).  \]
Thus, for $\Phi \in \mathbf{MVal}_{n-1}^{+,\infty}$,
representation (\ref{schu}) is equivalent to
\begin{equation} \label{schu1}
h(\Phi K,\cdot) = \mathrm{vol}_{n-1}(K|\,\cdot\,) \ast g_L,
\end{equation}
where $g_L \in C^{\infty}(S^{n-1})$ is the uniquely determined
even zonal function such that
\begin{equation} \label{schu2}
h(L,u)=\frac{1}{2}\int_{S^{n-1}}|u\cdot v|g_L(v)\,dv=(g_L\ast
\mbox{$\frac{1}{2}$}|\bar{e}\cdot\,.\,|)(u).
\end{equation}
Here, we have used that the convolution of zonal functions is
abelian.

\vspace{0.3cm}

\noindent {\bf Remark} Schneider \textbf{\cite{schneider74}},
Kiderlen \textbf{\cite{kiderlen05}}, and the author
\textbf{\cite{Schu06a}} originally \linebreak considered
translation invariant and $\mathrm{SO}(n)$ equivariant Minkowski
valuations of degree $1$ and $n-1$, respectively. However, their
results imply that for $n \geq 3$, these Minkowski valuations are
actually also $\mathrm{O}(n)$ equivariant.

\pagebreak

\centerline{\large{\bf{ \setcounter{abschnitt}{6}
\arabic{abschnitt}. Proof of the main results}}}

\vspace{0.7cm} \reseteqn \alpheqn \setcounter{koro}{0}

After these preparations, we are now in a position to give the
proof of Theorem \ref{thm1} and the stronger result Theorem
\ref{stronger}. At the end of this section we prove Theorem
\ref{thm2}.

The following result is a refined version of Theorem \ref{thm1}.

\begin{satz} \label{thm1rev} Suppose that $\Phi_i \in \mathbf{MVal}_i^{+,\infty}$, $1 \leq i
\leq n - 1$. Then there exists a smooth $\mathrm{O}(i) \times
\mathrm{O}(n - i)$ invariant measure $\mu$ on $S^{n-1}$ such that
for every $K \in \mathcal{K}^n$,
\[h(\Phi_iK,\cdot) = \mathrm{vol}_i(K|\,\cdot\,) \ast \mu.  \]
The measure $\mu$ can be chosen uniquely from a certain subset
$\mathbf{t}_i^{\infty}$ of $C^{\infty}(S^{n-1})$.
\end{satz}
{\it Proof.} Let $\varphi_i \in \mathbf{Val}$ denote the
associated real-valued valuation of $\Phi_i$. Since $\Phi_i \in
\mathbf{MVal}_i^{+,\infty}$, we have $\varphi_i \in
\mathbf{Val}_i^{+,\infty}$. Thus, it follows from Corollary
\ref{crofkoro} that there exists a unique smooth measure $\sigma
\in \mathbf{T}_i^{\infty}$ such that
\begin{equation} \label{crofton7}
\varphi_i(K)=\int_{\mathrm{Gr}_{i,n}}\mathrm{vol}_i(K|E)\,d\sigma(E),
\qquad K \in \mathcal{K}^n.
\end{equation}
By definition (\ref{real}) and the $\mathrm{O}(n)$ equivariance
of $\Phi_i$, we have for $K \in \mathcal{K}^n$ and every
$\vartheta \in \mathrm{O}(n - 1)$,
\[\varphi_i(\vartheta K)=h(\Phi_i(\vartheta K),\bar{e})=h(\vartheta \Phi_iK,\bar{e})=h(\Phi_iK,\bar{e})=\varphi_i(K).  \]
Therefore, the valuation $\varphi_i$ is $\mathrm{O}(n-1)$
invariant. By (\ref{crofton7}), this $\mathrm{O}(n-1)$ invariance
carries over to the measure $\sigma$.

We define the set $\mathbf{t}_i^{\infty} \subseteq
C^{\infty}(S^{n-1})$ by
\[\mathbf{t}_i^{\infty} = \{ \widehat{f}: f \in \mathbf{T}_i^{\infty}\, \mathrm{O}(n-1) \mbox{ invariant}\}.  \]
Since for any $\mathrm{O}(n-1)$ invariant $f \in
\mathbf{T}_i^{\infty}$, we can identify $\widehat{f}$ with an
$\mathrm{O}(n-1)$ right invariant and $\mathrm{O}(i) \times
\mathrm{O}(n - i)$ left invariant function on $\mathrm{O}(n)$,
the set $\mathbf{t}_i^{\infty}$ is well defined and consists of
$\mathrm{O}(i) \times \mathrm{O}(n - i)$ invariant functions in
$C^{\infty}(S^{n-1})$.

If we set $\mu = \widehat{\sigma}$, then, by the
$\mathrm{O}(n-1)$ invariance of $\sigma$, we have $\mu \in
\mathbf{t}_i^{\infty}$. Moreover, by (\ref{heta}) and
(\ref{crofton7}), it follows that for $\bar{\eta} \in S^{n-1}$
and every $K \in \mathcal{K}^n$,
\[h(\Phi_iK,\bar{\eta})= \varphi_i(\eta^{-1}K)=
\int_{\mathrm{Gr}_{i,n}}\mathrm{vol}_i(K|\eta
E)\,d\sigma(E)=(\mathrm{vol}_i(K|\cdot) \ast \mu)(\bar{\eta}),  \]
which concludes the proof of the theorem. \hfill $\blacksquare$

\pagebreak

We call a measure $\mu$ on $S^{n-1}$ a {\it Crofton measure} for
the Minkowski \linebreak valuation $\Phi_i \in \mathbf{MVal}_i^+$
if $h(\Phi_iK,\cdot) = \mathrm{vol}_i(K|\,\cdot\,) \ast \mu$. In
this case, we say $\Phi_i$ admits the Crofton measure $\mu$.

\vspace{0.3cm}

\noindent {\bf Examples:}
\begin{enumerate}
\item[(a)] For $i \in \{1, \ldots, n - 1\}$, we have $\Lambda_i \in
\mathbf{MVal}_i^{+,\infty}$. By (\ref{viwi}) and (\ref{wigrin}),
the Crofton measure of $\Lambda_i$ is a multiple of spherical
Lebesgue measure.
\item[(b)] For $i \in \{1, \ldots, n - 1\}$, it is well known (see, e.g.,
\textbf{\cite[\textnormal{p.\ 428}]{goodeyweil92}}) that the $i$th
projection operator $\Pi_i \in \mathbf{MVal}_i^+$ can be defined
by
\[h(\Pi_iK,\cdot) = \frac{\kappa_{n-1}}{\kappa_i} \mathrm{R}_{n-i}\mathrm{vol}_i^{\bot}(K|\,\cdot\,),  \]
where for $f \in C(\mathrm{Gr}_{i,n})$, the function $f^{\bot} \in
C(\mathrm{Gr}_{n-i,n})$ is defined by $f^{\bot}(E)=f(E^{\bot})$,
$E \in \mathrm{Gr}_{n-i,n}$. Here, $\mathrm{R}_i:
C(\mathrm{Gr}_{i,n}) \rightarrow C(S^{n-1})$ is the Radon
transform defined for $u \in S^{n-1}$ by
\[(\mathrm{R}_if)(u)=\int_{u \in E}f(E)\,dE, \qquad E \in \mathrm{Gr}_{i,n}.  \]
Let $\bar{E} \in \mathrm{Gr}_{i,n}$ denote the stabilizer of
$\mathrm{O}(i) \times \mathrm{O}(n-i)$. Grinberg and Zhang
\textbf{\cite[\textnormal{Lemma 3.2}]{grinbergzhang99}} have
shown that
\[\mathrm{R}_if=f \ast \mu_{S^{i-1}},  \]
where $\mu_{S^{i-1}}$ is the probability measure on $S^{n-1}$
uniformly concentrated on $S^{i-1}=S^{n-1} \cap \bar{E}$.
Consequently, we obtain
\[h(\Pi_iK,\cdot) = \frac{\kappa_{n-1}}{\kappa_i} \mathrm{vol}_i(K|\,\cdot\,) \ast \mu_{S^{n-i-1}}^{\bot},  \]
where
$\widehat{\mu^{\bot}}_{\hspace{-0.2cm}\!S^{n-i-1}}=\widehat{\mu}_{S^{n-i-1}}^{\bot}$.
\item[(c)] For $i \in \{2, \ldots, n - 1\}$, the $i$th mean
section operator $\mathrm{M}_i \in \mathbf{MVal}_{n+1-i}$ is not
even. However, Goodey and Weil \textbf{\cite[\textnormal{Theorem
5}]{goodeyweil92}} have shown that (for a suitable constant
$c_{n,i}$)
\[h(\mathrm{M}_iK,\cdot) + h(\mathrm{M}_i(-K),\cdot)=c_{n,i}\,\mathrm{R}_{n+1-i}\mathrm{vol}_{n+1-i}(K|\,\cdot\,).  \]
Thus, a multiple of $\mu_{S^{n-i}}$ is a Crofton measure for the
even part of $\mathrm{M}_i$.
\end{enumerate}

\vspace{0.3cm}

We note that not every Minkowski valuation $\Phi_i \in
\mathbf{MVal}_i^+$, $1 \leq i \leq n$, admits a Crofton measure.
For example, suppose that $\Phi \in \mathbf{MVal}_{n-1}^+$. Then,
by Theorem \ref{schuster}, there exists an origin-symmetric body
of revolution $L \in \mathcal{K}^n$ such that
\[h(\Phi K,\cdot) = S_{n-1}(K,\cdot) \ast h(L,\cdot).  \]
It follows from (\ref{schu1}) and (\ref{schu2}) that $\Phi$
admits a Crofton measure if and only if $L$ is a generalized
zonoid. However, it is well known that there exist convex bodies
(of revolution) which are not generalized zonoids.

\vspace{0.3cm}

An important additional property of Crofton measures of Minkowski
\linebreak valuations is contained in the following result:

\begin{satz} \label{support} Let $\Phi_i \in \mathbf{MVal}_i^{+}$, $1
\leq i \leq n - 1$. If $\mu$ is a Crofton measure for the
Minkowski valuation $\Phi_i$, then there exists an $\mathrm{O}(i)
\times \mathrm{O}(n-i)$ invariant convex body $L \in
\mathcal{K}^n$ such that
\[h(L,\cdot) = \widehat{\mathrm{C}_i \widehat{\mu}}.  \]
\end{satz}
{\it Proof.} Let $\pi_1: \mathrm{O}(n) \rightarrow
\mathrm{Gr}_{i,n}$ and $\pi_2: \mathrm{O}(n) \rightarrow S^{n-1}$
denote the canonical projections and let $e \in \mathrm{O}(n)$ be
the identity map. We denote by $\bar{E}=\pi_1(e)$ and
$\bar{e}=\pi_2(e)$ the stabilizers of $\mathrm{O}(i) \times
\mathrm{O}(n-i)$ and $\mathrm{O}(n-1)$, respectively.

Choose an $i$-dimensional subspace $F \in \mathrm{Gr}_{i,n}$ and
let $K \subseteq F$ be a convex body. For $u \in S^{n-1}$, a
theorem of Hadwiger \textbf{\cite[\textnormal{p.\
79}]{hadwiger51}} implies that
\[h(\Phi_iK,u)=f_i(F,u)\,\mathrm{vol}_i(K).  \]
This defines a continuous function $f_i: \mathrm{Gr}_{i,n} \times
S^{n-1} \rightarrow \mathbb{R}$ satisfying the following
properties:
\begin{enumerate}
\item[(a)] For each $F \in \mathrm{Gr}_{i,n}$, the function
$f_i(F,\cdot) \in C(S^{n-1})$ is the support function of a convex
body $L(F)$.
\item[(b)] The function $f_i(\cdot,\bar{e}) \in C(\mathrm{Gr}_{i,n})$
is the Klain function of the associated real-valued valuation
$\varphi_i$ of $\Phi_i$.
\item[(c)] For every $\vartheta \in \mathrm{O}(n)$, we have
$f_i(\vartheta F,u)=f_i(F,\vartheta^{-1}u)$.
\end{enumerate}
Define functions $g_1 \in C(\mathrm{Gr}_{i,n})$ and $g_2 \in
C(S^{n-1})$ by
\[g_1(F)=f_i(F,\bar{e}),\,\, F \in \mathrm{Gr}_{i,n} \qquad \mbox{and} \qquad g_2(u)=f_i(\bar{E},u),\,\, u \in S^{n-1}.   \]
From properties (a) and (c) of $f_i$, we deduce that $g_1$ is an
$\mathrm{O}(n-1)$ invariant function on $\mathrm{Gr}_{i,n}$ and
$g_2$ is an $\mathrm{O}(i) \times \mathrm{O}(n-i)$ invariant
support function of a convex body $L$. Moreover, property (c) of
$f_i$ implies that for every $\vartheta \in \mathrm{O}(n)$, we
have $g_1(\pi_1(\vartheta))=g_2(\pi_2(\vartheta^{-1}))$.
Therefore, we deduce that $g_1=\widehat{g_2}$.

From property (b) of $f_i$, (\ref{crofton7}) and (\ref{ciki}), we
finally obtain
\[h(L,\cdot)=h(\Phi_iK_{\bar{E}},\cdot)=\widehat{g_1}=\widehat{\mathrm{K}_i\varphi_i}=\widehat{\mathrm{C}_i\widehat{\mu}},  \]
where $K_{\bar{E}}$ is any convex body in $\bar{E}$ such that
$\mathrm{vol}_i(K_{\bar{E}})=1$. \hfill $\blacksquare$

\vspace{0.3cm}

Suppose that $\Phi_i \in \mathbf{MVal}_i^+$, $1 \leq i \leq n -
1$, admits a Crofton measure. \linebreak Then, by Theorem
\ref{support}, the Klain function of the  associated real-valued
\linebreak valuation of $\Phi_i$ is essentially the support
function of a convex body $L$. In particular, the convex body $L$
determines $\Phi_i$ uniquely.

\vspace{0.3cm}

Using Lemma \ref{lemdec}, we can prove a generalization of Theorem
\ref{thm1}:

\begin{satz} \label{stronger} Suppose that $\Phi \in \mathbf{MVal}^{+,\infty}$.
Then there exist constants $c_0, c_n \geq 0$ and smooth
$\mathrm{O}(i) \times \mathrm{O}(n - i)$ invariant measures
$\mu_i$ on $S^{n-1}$, where $1 \leq i \leq n - 1$, such that for
every $K \in \mathcal{K}^n$,
\[h(\Phi K,\cdot) = c_0 + \sum \limits_{i = 1}^{n-1} \mathrm{vol}_i(K|\,\cdot\,) \ast \mu_i + c_nV(K).  \]
The measures $\mu_i$ can be chosen uniquely from
$\mathbf{t}_i^{\infty}$.
\end{satz}
{\it Proof.} Let $\varphi \in \mathbf{Val}^{+,\infty}$ be the
associated real-valued valuation of $\Phi$. By (\ref{decsmooth}),
there exist smooth valuations $\varphi_i \in
\mathbf{Val}_i^{+,\infty}$, $0 \leq i \leq n$, such that
\linebreak $\varphi = \varphi_0 + \ldots + \varphi_n$. By Lemma
\ref{lemdec} on the other hand, there exist constants $c_0, c_n
\geq 0$ and even $g_i(K,\cdot) \in C(S^{n-1})$ such that
\[h(\Phi K,\cdot) = c_0 + \sum \limits_{i = 1}^{n-1} g_i(K,\cdot) + c_nV(K).  \]
Clearly, for $i \in \{1, \ldots, n - 1\}$, we have
$\varphi_i(K)=g_i(K,\bar{e})$. By Lemma \ref{lemdec},
$g_i(\vartheta K,u)=g_i(K,\vartheta^{-1}u)$ for every $\vartheta
\in \mathrm{O}(n)$. Thus, it follows, as in the proof of Theorem
\ref{thm1rev}, that there exist unique measures $\mu_i \in
\mathbf{t}_i^{\infty}$ such that
\[g_i(K,\cdot) =  \mathrm{vol}_i(K|\,\cdot\,) \ast \mu_i. \]
\vspace{-1cm}

\hfill $\blacksquare$

\vspace{0.3cm}

We note that the techniques applied in the proof of Theorem
\ref{thm1rev} together with a straightforward adaptation of Lemma
\ref{lemdec} immediately yield (the proof is almost verbatim the
same as the proof of Theorem \ref{stronger}) a characterization
of valuations with values in the space of continuous functions on
$S^{n-1}$:

\begin{koro} A map $f: \mathcal{K}^n \rightarrow C(S^{n-1})$ is a
smooth translation invariant \linebreak and $\mathrm{O}(n)$
equivariant even valuation if and only if there exist constants
$c_0, c_n \in \mathbb{R}$ and smooth $\mathrm{O}(i) \times
\mathrm{O}(n - i)$ invariant measures $\mu_i$ on $S^{n-1}$, where
$1 \leq i \leq n - 1$, such that for every $K \in \mathcal{K}^n$,
\[f(K,\cdot) = c_0 + \sum \limits_{i = 1}^{n-1} \mathrm{vol}_i(K|\,\cdot\,) \ast \mu_i + c_nV(K).  \]
\end{koro}

Here, the continuous translation invariant and $\mathrm{O}(n)$
equivariant valuation $f: \mathcal{K}^n \rightarrow C(S^{n-1})$
is called smooth, if the associated real-valued valuation $\psi
\in \mathbf{Val}$, defined by $\psi(K)=f(K,\bar{e}), K \in
\mathcal{K}^n$, is smooth.

\vspace{0.3cm}

We conclude this section with the proof of Theorem \ref{thm2}:

\begin{satz} \label{thm2rev} For every Minkowski valuation $\Phi \in \mathbf{MVal}^+$, there exists a sequence $\Phi^m \in
\mathbf{MVal}^{+,\infty}$, $m \in \mathbb{N}$, such that $\Phi^m$
converges to $\Phi$ uniformly on compact subsets of
$\mathcal{K}^n$.
\end{satz}
{\it Proof.} Let $f_m \in C^{\infty}(S^{n-1})$ be a spherical
approximate identity. For each $m \in \mathbb{N}$, we define a
continuous map $\Phi^m: \mathcal{K}^n \rightarrow \mathcal{K}^n$
by
\[h(\Phi^mK,\cdot) = h(\Phi K,\cdot) \ast f_m, \qquad K \in \mathcal{K}^n.  \]
It follows from \textbf{\cite[\textnormal{Proposition
3.2}]{kiderlen05}} that the spherical convolution from the right
with nonnegative functions and measures maps support functions to
support functions. Therefore, $\Phi^m$ is well defined. Moreover,
it is easy to verify that $\Phi^m$ is an even translation
invariant Minkowski valuation. By (\ref{roteq}), $\Phi^m$ is also
$\mathrm{O}(n)$ equivariant. Consequently, $\Phi^m \in
\mathbf{MVal}^+$.

Let $\varphi^m \in \mathbf{Val}^{+}$ be the associated real-valued
valuation of $\Phi^m$ and let
\[\varphi^m=c_0^m + \varphi^m_1 + \ldots + \varphi_{n-1}^m+c_n^mV,\]
where $c_0^m, c_n^m \in \mathbb{R}$ and $\varphi^m_i \in
\mathbf{Val}^+_i$, $1 \leq i \leq n - 1$, be the decomposition of
$\varphi^m$ into homogeneous parts (which follows from Theorem
\ref{mcmullen}). In order to show that for each $m \in
\mathbb{N}$, the Minkowski valuation $\Phi^m$ is smooth, it
suffices to show that $\varphi^m_i$, $1 \leq i \leq n - 1$, is
smooth.

\pagebreak

Since $f_m \in C^{\infty}(S^{n-1})$, we deduce from an application
of Lemma \ref{lemdec} to $\Phi$ and $\Phi^m$, that for each $m \in
\mathbb{N}$ and every $K \in \mathcal{K}^n$, there exist even
$g^m_i(K,\cdot) \in C^{\infty}(S^{n-1})$ such that
\[h(\Phi^m K,\cdot) = c_0^m + \sum \limits_{i = 1}^{n-1} g^m_i(K,\cdot) + c_n^mV(K).  \]
Choose an $i$-dimensional subspace $F \in \mathrm{Gr}_{i,n}$ and
let $K \subseteq F$ be a convex body. As in the proof of Theorem
\ref{support}, it follows that for each $m \in \mathbb{N}$, there
exists a continuous function $\zeta^m_i: \mathrm{Gr}_{i,n} \times
S^{n-1} \rightarrow \mathbb{R}$ such that
\[g_i^m(K,u)=\zeta^m_i(F,u)\,\mathrm{vol}_i(K|F).  \]
In fact, since $g_i^m(K,\cdot) \in C^{\infty}(S^{n-1})$, we also
have $\zeta^m_i(F,\cdot) \in C^{\infty}(S^{n-1})$ for each $F \in
\mathrm{Gr}_{i,n}$. Therefore, it follows from the proof of
Theorem \ref{support} that
\[\mathrm{K}_i \varphi_i^m = \widehat{\zeta^m(\bar{E},\cdot)}  \]
is smooth. Consequently, $\varphi_i^m$, $1 \leq i \leq n - 1$, is
smooth which in turn implies that $\Phi^m \in
\mathbf{MVal}^{+,\infty}$.

It remains to show that $\Phi^m$ converges to $\Phi$ uniformly on
compact subsets. By (\ref{zonsuff}), we may assume that $f_m$ is
zonal for each $m \in \mathbb{N}$. If $g_i(K,\cdot) \in
C(S^{n-1})$ denotes the degree $i$ component of $h(\Phi K,\cdot)$,
then $g_i^m(K,\cdot) = g_i(K,\cdot) \ast f_m$. By Lemma
\ref{approxid}, $g_i^m(K,\cdot)$ converges uniformly to
$g_i(K,\cdot)$ for each $K \in \mathcal{K}^n$. Moreover, it is
not hard to show that
\[|\varphi_i^m(K) - \varphi_i(K)| = \left | \int_{S^{n-1}} g_i(K,u)f_m(u)\,du - g_i(K,\bar{e}) \right |,  \]
where $\varphi_i$ denotes the degree $i$ part of the associated
real-valued valuation $\varphi$ of $\Phi$. Hence, the
$\varphi_i^m$ converge to $\varphi_i$ pointwise. Since the map $K
\mapsto g_i(K,\cdot)$ is uniformly continuous on every compact
subset of $\mathcal{K}^n$ and
\[|\varphi_i^m(K) - \varphi_i^m(L)| \leq \|g_i(K,\cdot)-g_i(L,\cdot)\|_{\infty},  \]
the $\varphi_i^m$ are equicontinuous on every compact subset of
$\mathcal{K}^n$. It follows that the $\varphi_i^m$ converge to
$\varphi_i$ and thus $\varphi^m$ to $\varphi$ uniformly on
compact subsets of $\mathcal{K}^n$. Since $\mathrm{O}(n)$ is
compact, it is easy to verify, using (\ref{heta}), that this
implies uniform convergence of $\Phi^m$ to $\Phi$ on compact
subsets of $\mathcal{K}^n$. \hfill $\blacksquare$

\pagebreak

\centerline{\large{\bf{ \setcounter{abschnitt}{7}
\arabic{abschnitt}. A Brunn--Minkowski type inequality}}}

\vspace{0.7cm} \reseteqn \alpheqn \setcounter{koro}{0}

As an application of Theorem \ref{thm1}, we present in this last
section the proof of Theorem \ref{thm3}. It is based on the
techniques developed by Lutwak in \textbf{\cite{lutwak93}}.

From (\ref{steiner}) and Theorem 5.2, the author deduced in
\textbf{\cite{Schu06a}} that for any $\Psi \in
\mathbf{MVal}_{n-1}$, there exist {\it derived} Minkowski
valuations $\Psi_i \in \mathbf{MVal}_i$, where $0 \leq i \leq n -
1$, such that for every $K \in \mathcal{K}^n$,
\[\Psi(K + \varepsilon B)= \sum \limits_{i=0}^{n-1} \varepsilon^{n-1-i} {n - 1 \choose i} \Psi_i K.  \]
Moreover, the author obtained in \textbf{\cite{Schu06b}} an array
of geometric inequalities for the intrinsic volumes of derived
(non-trivial) Minkowski valuations $\Psi_i$. In particular, the
following Brunn--Minkowski type inequality was established (cf.\
\textbf{\cite[\textnormal{Theorem 6.8}]{Schu06b}}): If $K, L \in
\mathcal{K}^n_\mathrm{o}$ and $2 \leq i \leq n-1$, $1 \leq j \leq
n$, then
\begin{equation} \label{genbmbmind}
V_j(\Psi_i(K+L))^{1/ij} \geq
V_j(\Psi_iK)^{1/ij}+V_j(\Psi_iL)^{1/ij},
\end{equation}
with equality if and only if $K$ and $L$ are homothetic.

We believe that inequality (\ref{genbmbmind}) holds in fact for
all Minkowski valuations in $\mathbf{MVal}_i$. Theorem \ref{thm3}
confirms this conjecture in the case of even valuations and $j =
i + 1$.

\vspace{0.2cm}

From now on we always assume that $\Phi_i \in \mathbf{MVal}_i^+$.
A critical ingredient in the proof of Theorem \ref{thm3} is the
following Lemma.

\begin{lem} \label{lemdurch1} If $K, L \in \mathcal{K}^n$ and $1 \leq i \leq
n-1$, then
\begin{equation} \label{durch1a}
W_{n-1-i}(K,\Phi_iL)=W_{n-1-i}(L,\Phi_iK).
\end{equation}
In particular, there exists a constant $r(\Phi_i) \geq 0$ such
that
\begin{equation} \label{durch1b}
W_{n-1}(\Phi_i K)=r(\Phi_i)W_{n-i }(K).
\end{equation}
\end{lem}
{\it Proof:} By Theorem \ref{thm2}, it suffices to prove the
statement for $\Phi_i \in \mathbf{MVal}_i^{+,\infty}$. For $K \in
\mathcal{K}^n$ and $1 \leq i \leq n - 1$, we define a measure
$s_i(K,\cdot)$ on $S^{n-1}$ by
\[s_i(K,\cdot)=\frac{1}{2}S_i(K,\cdot) + \frac{1}{2}S_i(-K,\cdot).  \]
Since $\Phi_i$ is even, it follows from (\ref{wisi}) that
\begin{equation} \label{phisi}
W_{n-1-i}(K,\Phi_iL)=\frac{1}{n}\int_{S^{n-1}}h(\Phi_iL,u)\,ds_i(K,u).
\end{equation}

Let $f_m \in C^{\infty}(S^{n-1})$ be a spherical approximate
identity. For each $m \in \mathbb{N}$, we define a function
$s_i^m(K,\cdot) \in C^{\infty}(S^{n-1})$ by
\[s_i^m(K,\cdot) = s_i(K,\cdot) \ast f_m.  \]
From Lemma \ref{approxid}, (\ref{roteq}) and well-known properties
of the area measures $S_i(K,\cdot)$ (see, e.g.,
\textbf{\cite[\textnormal{Chapter 4 \& 5}]{schneider93}}), it
follows that for each $m \in \mathbb{N}$:
\begin{enumerate}
\item[(a)] The map $K \mapsto s_i^m(K,\cdot)$ is a continuous,
translation invariant even valuation of degree $i$.
\item[(b)] For every $\vartheta \in \mathrm{O}(n)$ and $K \in
\mathcal{K}^n$, we have $s_i^m(\vartheta
K,u)=s_i^m(K,\vartheta^{-1}u)$.
\end{enumerate}
Moreover, as in the proof of Theorem \ref{thm2rev}, one can show
that for each $m \in \mathbb{N}$, the real-valued valuation $K
\mapsto s_i^m(K,\bar{e})$ is smooth. Thus, it follows, as in the
proof of Theorem \ref{thm1rev}, that there exist unique measures
$\sigma_i^m \in \mathbf{t}_i^{\infty}$ such that
\[s_i^m(K,\cdot) =  \mathrm{vol}_i(K|\,\cdot\,) \ast \sigma_i^m. \]
Thus, from (\ref{phisi}), (\ref{wisi}), Lemma \ref{approxid} and
Theorem \ref{thm1}, we obtain
\begin{equation} \label{form1}
W_{n-1-i}(K,\Phi_iL)=\lim_{m} \frac{1}{n} \left \langle
\mathrm{vol}_i(L|\,\cdot\,) \ast \mu_i ,
\mathrm{vol}_i(K|\,\cdot\,) \ast \sigma_i^m \right \rangle,
\end{equation}
for some $\mu_i \in \mathbf{t}_i^{\infty}$. From (\ref{adjoint})
and (\ref{commute}), it follows that
\begin{equation} \label{form2}
\left \langle \mathrm{vol}_i(L|\,\cdot\,) \ast \mu_i ,
\mathrm{vol}_i(K|\,\cdot\,) \ast \sigma_i^m \right \rangle=\left
\langle \mathrm{vol}_i(L|\,\cdot\,), \mathrm{vol}_i(K|\,\cdot\,)
\ast \widehat{\mu_i \ast \widehat{\sigma_i^m}}  \right \rangle.
\end{equation}
Since the measure $\mu_i$ is left $\mathrm{O}(i) \times
\mathrm{O}(n-i)$ invariant and the measure $\widehat{\sigma_i^m}$
is right $\mathrm{O}(i) \times \mathrm{O}(n-i)$ invariant, we
obtain from Lemma \ref{selfad}, that
\[\widehat{\mu_i \ast \widehat{\sigma_i^m}} = \mu_i \ast \widehat{\sigma_i^m}  \]
Combining, (\ref{form1}), (\ref{form2}) and (\ref{adjoint}), we
obtain
\[W_{n-1-i}(K,\Phi_iL)=\lim_{m} \frac{1}{n}
\left \langle \mathrm{vol}_i(L|\,\cdot\,) \ast \sigma_i^m,
\mathrm{vol}_i(K|\,\cdot\,) \ast \mu_i \right
\rangle=W_{n-1-i}(L,\Phi_iK).\]

Finally, note that, by the $\mathrm{O}(n)$ equivariance of
$\Phi_i$, we have $\Phi_iB=r(\Phi_i)B$, for some constant
$r(\Phi_i) \geq 0$. Thus,
\[W_{n-1}(\Phi_iK)=W_{n-1-i}(B,\Phi_iK)=W_{n-1-i}(K,\Phi_iB)=r(\Phi_i)W_{n-i}(K).  \]
\hfill $\blacksquare$ \vspace{0.3cm}

By (\ref{viwi}), our next result is equivalent to Theorem
\ref{thm3}:

\begin{satz} \label{satzgenbmbm} If $K, L \in \mathcal{K}^n_{\mathrm{o}}$ and $2 \leq i \leq n-1$, then
\begin{equation*}
W_{n-1-i}(\Phi_i(K+L))^{1/i(i+1)} \geq
W_{n-1-i}(\Phi_iK)^{1/i(i+1)}+W_{n-1-i}(\Phi_iL)^{1/i(i+1)}.
\end{equation*}
If $\Phi_i \mathcal{K}^n_{\mathrm{o}} \subseteq
\mathcal{K}^n_{\mathrm{o}}$, then equality holds if and only if
$K$ and $L$ are homothetic.
\end{satz}
{\it Proof:} By (\ref{durch1a}) and (\ref{mostgenbm}), we have for
$Q \in \mathcal{K}^n_{\mathrm{o}}$,
\begin{eqnarray*}
W_{n-1-i}(Q,\Phi_i(K+L))^{1/i} & = & W_{n-1-i}(K+L,\Phi_iQ)^{1/i}
\\
& \geq & W_{n-1-i}(K,\Phi_i Q)^{1/i}+W_{n-1-i}(L,\Phi_i Q)^{1/i}
\\
&=& W_{n-1-i}(Q,\Phi_i K)^{1/i}+W_{n-1-i}(Q,\Phi_i L)^{1/i}.
\end{eqnarray*}
From inequality (\ref{genmink}), we further deduce that
\begin{equation} \label{mink1}
W_{n-1-i}(Q,\Phi_iK)^{i+1}\geq
W_{n-1-i}(Q)^{i}W_{n-1-i}(\Phi_iK),
\end{equation}
and
\begin{equation} \label{mink2}
W_{n-1-i}(Q,\Phi_iL)^{i+1}\geq W_{n-1-i}(Q)^{i}W_{n-1-i}(\Phi_iL).
\end{equation}
Thus, if we set $Q=\Phi_i(K+L)$, we obtain the desired inequality
\begin{equation*}
W_{n-1-i}(\Phi_i(K+L))^{1/i(i+1)} \geq
W_{n-1-i}(\Phi_iK)^{1/i(i+1)}+W_{n-1-i}(\Phi_iL)^{1/i(i+1)}.
\end{equation*}
Suppose now that equality holds and that $\Phi_i
\mathcal{K}^n_{\mathrm{o}} \subseteq \mathcal{K}^n_{\mathrm{o}}$.
Since $\Phi_i K$ is origin-symmetric for every $K \in
\mathcal{K}^n$, we deduce from the equality conditions of
(\ref{mink1}) and (\ref{mink2}), that there exist $\lambda_1,
\lambda_2 > 0$ such that
\begin{equation} \label{hom2}
\Phi_iK=\lambda_1\Phi_i(K+L) \qquad \mbox{and} \qquad
\Phi_iL=\lambda_2 \Phi_i(K+L)
\end{equation}
and
\[\lambda_1^{1/i}+\lambda_2^{1/i} = 1.   \]
Moreover, since $r(\Phi_i) > 0$, (\ref{durch1b}) and (\ref{hom2})
imply
\[W_{n-i}(K)=\lambda_1W_{n-i}(K+L) \qquad \mbox{and} \qquad  W_{n-i}(L)=\lambda_2W_{n-i}(K+L).  \]
Hence, we have
\[W_{n-i}(K+L)^{1/i}=W_{n-i}(K)^{1/i}+W_{n-i}(L)^{1/i},  \]
which implies, by (\ref{quermassbm}), that $K$ and $L$ are
homothetic.
\hfill $\blacksquare$\\

\medskip\noindent{\bf Acknowledgements.}
The work of the author was supported by the Austrian Science Fund
(FWF), within the project ``Valuations on Convex Bodies", Project
Number: P\,18308.

\vspace{0.5cm}

Vienna University of Technology \par Institute of Discrete
Mathematics and Geometry \par Wiedner Hauptstra\ss e 8--10/1046
\par A--1040 Vienna, Austria \par
franz.schuster@tuwien.ac.at

\end{document}